\documentclass[11pt]{article}
\usepackage{dsfont, amssymb,amsmath,amscd,latexsym, amsthm, amsxtra,amsfonts, mathrsfs}
\usepackage[all]{xy}
\begin{document}
\baselineskip 20pt
\title{{ {  Constituting Atoms of a $\sigma$ Algebra via Its
Generator}}}
\author{Jinshan \ Zhang\thanks{Electronic address: zjs02@mails.tsinghua.edu.cn}
\\ \small Department of Mathematical Sciences,
\  Tsinghua University,\\
\small Beijing 100084,  China }
\date{}
\maketitle
\vskip10pt
\begin{center} {\bf Abstract }\end{center}
\ \ \;In this paper, a very weak sufficient condition for
determining atoms by the generator is presented. The condition,
though not being a necessary one, is shown to be almost the weakest
one in the sense that it can hardly be improved.
  \text{}
\section{ Introduction }

\ \ \;The atom is an advanced concept in modern probability theory.
An intensive application of the properties of atoms plays an
important role in one of the proofs of the famous  Vatali-Hahn-Saks
Theorem \cite{Yan04}. This concept is essential  in the development
of conditional probability during the recent
decades\cite{BD75,BR63,Haj03,SSK01}. However when atom is used, the
relations between $\sigma$ algebra and its atoms is mostly assumed
to be known. Namely most of the authors consider such relations as
conditions in their research. This paper focuses on how to determine
such relations, which means to find the structures and  the
constructions of atoms of a $\sigma$ algebra. The most relevant
result to this paper should be Blackwell theorem \cite{DM78}, which
is quite useful and has reduced inclusion between Blackwell $\sigma$
algebras to comparing their atoms. Hence, determining the atoms of a
$\sigma$ algebra becomes so significant when applying that powerful
theorem to solve problems.

The $\sigma$ algebra itself, however, is usually too large or can
not be efficiently obtained by  given its generator only. To solve
this matter, reducing constituting atoms from the $\sigma$ algebra
to constituting them from the generator would be a feasible and
efficient approach. Therefore, at least some sufficient conditions
for the generator on its structure, which can be used to constitute
the atoms of the $\sigma$ algebra generated, should be provided. The
condition for the generator proposed in this paper only requires
$\kappa(\mathcal{C})=\sigma(\mathcal{C})$. This condition, though
not being a necessary one, is shown to be almost the weakest one in
the sense that it can hardly be improved.

The following is a brief road-map of the paper. In section 2, some
preliminary knowledge is introduced, which may be used in this
paper. We go on in section 3 to prove the main theorem of the paper.
Some useful corollaries following the theorem construct the subject
of section 4. In section 5, the condition in the main theorem is
discussed through specific examples and theoretical analysis.
Concluding remarks are proposed.

\section{Preliminaries}

\ \ \;In this section we review some fundamental  concepts and
related results that would  be utilized in this paper. They are
mostly important definitions and theorems in probability theory. In
the following, we always let $\Omega$ be a set,  and $\mathcal{C}$
be a collection (set) of subsets of the set $\Omega$.\\\\
\textbf{Definition 2.1} $\mathcal{C}$ is called a monotone class, if
$A_{n}\in\mathcal{C},n\geq 1, A_{n}\uparrow A$ or $A_{n}\downarrow
A$ $\Rightarrow$ $A\in\mathcal{C}$. Further, $\mathcal{C}$ is called
a $\lambda$ class, if  $\mathcal{C}$ is a monotone class,
$\Omega\in\mathcal{C}$, and $\forall A$,$B\in\mathcal{C}$, $B\subset
A$ $\Rightarrow$ $A\cap B^{c}\in \mathcal{C}$.\\\\
\textbf{Definition 2.2} Let $(\Omega, \mathcal{F})$ be a measurable
space, where $\mathcal{F}$ is a $\sigma$ algebra on $\Omega$. For
any $\omega\in\Omega$ define
$\mathcal{F}_{\omega}=\{B\in\mathcal{F}| \omega\in B\}$. Then
$A(\omega)=\bigcap\limits_{B\in\mathcal{F}_{\omega}}B$ is called an
atom of $\mathcal{F}$ containing $\omega$. \\\\
Let $\mathcal{C}_{\cap f}=\{A|$ $A=\bigcap\limits_{i=1}^{n}A_{i} ,
A_{i}\in\mathcal{C}, i=1,\cdots,n,n\geq1\}$ be the set closed under
finite intersection. Similarly, $\mathcal{C}_{\cup f}$,
$\mathcal{C}_{\sum f}$, $\mathcal{C}_{\delta}$,
$\mathcal{C}_{\sigma}$, $\mathcal{C}_{\sum \sigma}$ denote the sets
closed under finite union, finite disjoint union, countable
intersection, countable union, countable disjoint union
respectively.
$\mathcal{C}_{\sigma\delta}=(\mathcal{C}_{\sigma})_{\delta}$.
$\sigma(\mathcal{C})$, $\lambda(\mathcal{C})$, $m(\mathcal{C})$
denotes the minimum $\sigma$ algebra, $\lambda$ class,  and monotone
class containing $\mathcal{C}$ respectively and $\mathcal{C}$ is the
generator of them.

\section{Main theorem}

Let $\mathcal{C}$ be a collection of sets on $\Omega$ and
$\mathcal{F}=\sigma(\mathcal{C})$. For any $ \omega\in\Omega$,
define $\mathcal{C}_{\omega}=\{B\in \mathcal{C}|$
 $\omega\in B\}$ and $A(\omega)$ the atom containing $\omega$ of
$\mathcal{F}$. Question is: under what condition on $\mathcal{C}$,
there would be
$$A(\omega)=\bigcap\limits_{B\in\mathcal{C}_{\omega}}B. $$ Yan
\cite{Yan04} shows that this is true, if $\mathcal{C}$ is an
algebra. A much weaker condition on $\mathcal{C}$ is proposed in
this section, which is the main result of this paper. In order to
show the main result of this paper, we introduce the following
concept, which is created in this paper to show how good our
condition is.
\\\\
\textbf{Definition 3.1} $\mathcal{C}$ is called a $\kappa$ class, if
it is closed under countable intersection and countable union.\\\\
Denote by $\kappa(\mathcal{C})$ the minimum $\kappa$ class
containing $\mathcal{C}$, and $\mathcal{C}$ is called generator of
$\kappa(\mathcal{C})$. To complete the proof of the main
 result, we need the following lemma.\\
\textbf{Lemma 3.1} Let $\mathcal{C}$ be a collection of sets on
$\Omega$, $\mathcal{F}=\sigma(\mathcal{C})$. $\forall
\omega\in\Omega$, define $\mathcal{C}_{\omega}=\{B\in\mathcal{C}|$
$\omega\in B\}$ and
$A_{\mathcal{C}}(\omega)=\bigcap\limits_{B\in\mathcal{C}_{\omega}}B$.
Let $\mathcal{G}=\{B\in\mathcal{F}|$ $\omega\notin B$, or $\omega\in
B$ and $A_{\mathcal{C}}(\omega)\subset B\}$ $=$ $\{B\in\mathcal{F}|$
$\omega\notin B\}$ $\bigcup$ $\{B\in\mathcal{F}|$ $\omega\in B$,
$A_{\mathcal{C}}(\omega)\subset B\}. $\\
Then $\mathcal{G}$ satisfies the following three properties:\\
1). $\mathcal{C}\subset\mathcal{G}$,
$\mathcal{C}_{\sigma}\subset\mathcal{G}$.\\
2). $\mathcal{G}$ is closed under the operation of countable
union and countable intersection.\\
3). $\mathcal{G}$ is a $\kappa$ class. In particular, it is a
monotone
class.\\
\textbf{Proof:} Let
$\mathcal{C}_{\sigma\omega}=\{B\in\mathcal{C}_{\sigma}|$ $\omega\in
B\}$ and
$A_{\mathcal{C}_{\sigma}}(\omega)=\bigcap\limits_{B\in\mathcal{C}_{\sigma\omega}}B$.
Claim $\forall \omega\in\Omega$,
$A_{\mathcal{C}}(\omega)=A_{\mathcal{C}_{\sigma}}(\omega)$. First,
$\mathcal{C}_{\omega}\subset\mathcal{C}_{\sigma\omega}$, then
$A_{\mathcal{C}_{\sigma}}(\omega)=\bigcap\limits_{B\in\mathcal{C}_{\sigma\omega}}B$
$\subset$
$\bigcap\limits_{B\in\mathcal{C}_{\omega}}B=A_{\mathcal{C}}(\omega)$.
Consider the definition of $\mathcal{C}_{\sigma}$ and
$\mathcal{C}_{\sigma\omega}$, $\forall$
$B\in\mathcal{C}_{\sigma\omega}$, $\exists$
$\{A_{n}\}_{n=1}^{\infty}\subset\mathcal{C}$ such that
$B=\bigcup\limits_{n=1}^{\infty}A_{n}$. Hence, there exits $N$ such
that $\omega\in A_{N}$, then there are
$A_{\mathcal{C}}(\omega)\subset A_{N}\subset B$ and
$A_{\mathcal{C}}(\omega)\subset\bigcap\limits_{B\in\mathcal{C}_{\sigma\omega}}B=A_{\mathcal{C}_{\sigma}}(\omega)$.
Hence $A_{\mathcal{C}}(\omega)=A_{\mathcal{C}_{\sigma}}(\omega)$.
Now, let's prove the lemma. \\
For the property 1. $\forall$ $B\in\mathcal{C}$, if $\omega\notin B$
then $B\in \mathcal{G}$; Otherwise, if $\omega\in B$, since
$A_{\mathcal{C}}(\omega)=\bigcap\limits_{B\in\mathcal{C}_{\omega}}B$,
we have $A_{\mathcal{C}}(\omega)\subset B$, then $B\in\mathcal{G}$.
Hence $\mathcal{C}\subset\mathcal{G}$. Frow the claim, we know
$\mathcal{G}=\{B\in\mathcal{F}|$ $\omega\notin B$, or $\omega\in B$
and  $A_{\mathcal{C}_{\sigma}}(\omega)\subset B\}$ $=$
$\{B\in\mathcal{F}|$ $\omega\notin B\}$ $\bigcup$
$\{B\in\mathcal{F}|$ $\omega\in B$,
$A_{\mathcal{C}_{\sigma}}(\omega)\subset B\}$. Thus, similarly,
$\mathcal{C}_{\sigma}\subset\mathcal{G}$.\\
For the property 2. Suppose
$\{A_{n}\}_{n=1}^{\infty}\subset\mathcal{G}$.\\
(i). If $\forall n$, $\omega\notin A_{n}$, then
$\omega\notin\bigcup\limits_{n=1}^{\infty}A_{n}$. Hence
$\bigcup\limits_{n=1}^{\infty}A_{n}\in \mathcal{G}$.\\
(ii). If $\exists n$ such that $\omega\in A_{n}$ then
$A_{\mathcal{C}}(\omega)\subset
A_{n}\subset\bigcup\limits_{n=1}^{\infty}A_{n}$ . Obviously,
$\omega\in \bigcup\limits_{n=1}^{\infty}A_{n}$. Hence
$\bigcup\limits_{n=1}^{\infty}A_{n}\in\mathcal{G}$\\
Considering (i) and (ii), $\mathcal{G}$ is closed under countable
union.\\
(iii). If $\exists n$ such that $\omega\notin A_{n}$, then
$\omega\notin\bigcap\limits_{n=1}^{\infty}A_{n}$. Thus
$\bigcap\limits_{n=1}^{\infty}A_{n}\in\mathcal{G}$. \\
(iv). If $\forall n$, $\omega\in A_{n}$, then $\omega\in
\bigcap\limits_{n=1}^{\infty}A_{n}$. Since
$A_{\mathcal{C}}(\omega)\subset A_{n}$($\forall n$),
$A_{\mathcal{C}}(\omega)\subset \bigcap\limits_{n=1}^{\infty}A_{n}$.
Hence $\bigcap\limits_{n=1}^{\infty}A_{n}\in\mathcal{G}$.\\
Considering (iii) and (iv), $\mathcal{G}$ is closed under countable
intersection.\\
For property 3. From property 2, we know $\mathcal{G}$ is a $\kappa$
class. In particular, if $A_{n}\uparrow A$ then
$A=\bigcup\limits_{n=1}^{\infty}A_{n}$, and if $A_{n}\downarrow A$,
then $A=\bigcap\limits_{n=1}^{\infty}A_{n}$. Hence $\mathcal{G}$ is
a monotone class. \ \
$\Box$\\

Using Lemma 3.1, now we prove the main result of this paper.\\
\textbf{Theorem 3.1}  Let $\mathcal{C}$ be a collection of sets on
$\Omega$, $\mathcal{F}=\sigma(\mathcal{C})$ and
$A_{\mathcal{F}}(\omega)$ the atom of $\mathcal{F}$ containing
$\omega$. $\forall \omega\in\Omega$, define
$\mathcal{C}_{\omega}=\{B\in\mathcal{C}|$ $\omega\in B\}$ and
$A_{\mathcal{C}}(\omega)=\bigcap\limits_{B\in\mathcal{C}_{\omega}}B$.
If the generator $\mathcal{C}$ satisfies the property that $\forall$
$A\in\mathcal{C}$ $\Rightarrow$ $A^{c}\in \kappa(\mathcal{C})$, then
$$A_{\mathcal{F}}(\omega)=A_{\mathcal{C}}(\omega). $$
\textbf{Proof:} $\forall\omega$, let
$\mathcal{G}_{1}=\{B\in\mathcal{F}|$ $\omega\notin B$, or $\omega\in
B$ and $A_{\mathcal{C}}(\omega)\subset B\}$ and
$\mathcal{G}_{2}=\{A\in\mathcal{G}_{1}|$
$A^{c}\in\mathcal{G}_{1}\}$. Then $\mathcal{G}_{2}$ satisfies the following properties.\\
(a). $\forall$ $A\in\mathcal{C}$, $A^{c}\in \kappa(\mathcal{C})$, by
the property 1 and 3 of $\mathcal{G}_{1}$ in Lemma 3.1, we know
$\kappa(\mathcal{C})\subset\mathcal{G}_{1}$, then
$A^{c}\in\mathcal{G}_{1}$. Hence
$\mathcal{C}\subset\mathcal{G}_{2}$.\\
(b). Since $\mathcal{G}_{1}$ is a monotone class, it's easy to
check $\mathcal{G}_{2}$ is a monotone class. \\
(c). Now let's check $\mathcal{G}_{2}$ is an algebra.

(i). $\forall$ $A\in\mathcal{G}_{2}$, then $A\in\mathcal{G}_{1}$,
$A^{c}\in\mathcal{G}_{1}$, $(A^{c})^{c}\in\mathcal{G}_{1}$, hence
$A^{c}\in\mathcal{G}_{2}$.

(ii).$\forall$ $A$,$B\in\mathcal{G}_{2}$, then
$A$,$A^{c}\in\mathcal{G}_{1}$ and $B$,$B^{c}\in\mathcal{G}_{1}$.
Consider the property 2 of $\mathcal{G}_{1}$, we know $A\cap
 B\in\mathcal{G}_{1}$, $A^{c}\cup B^{c}\in \mathcal{G}_{1}$, then $(A\cap B)^{c}\in
 \mathcal{G}_{1}$. Hence $A\cap B\in\mathcal{G}_{2}$.\\
 Considering (i) and (ii), we show $\mathcal{G}_{2}$ is an algebra.
Now from (a), (b) and (c), $\mathcal{G}_{2}$ is a monotone class and
algebra containing $\mathcal{C}$. By Monotone Class Theorem,
$\mathcal{F}=\sigma(\mathcal{C})\subset\mathcal{G}_{2}$. Then
$\mathcal{G}_{2}=\mathcal{F}$.
$\mathcal{G}_{2}\subset\mathcal{G}_{1}\subset\mathcal{F}$, then
$\mathcal{G}_{1}=\mathcal{F}$. Noting that
$\mathcal{G}_{1}/\mathcal{F}_{\omega}=\{B\in\mathcal{F}|$
$\omega\notin B\}$(recall $\mathcal{F}_{\omega}=\{B\in\mathcal{F}|$
$\omega\in B\}$), then $\forall$$B\in\mathcal{F}_{\omega}$
$A_{\mathcal{C}}(\omega)\subset B$. Hence
$A_{\mathcal{C}}(\omega)\subset A_{\mathcal{F}}(\omega)$. since
$\mathcal{C}_{\omega}\subset\mathcal{F}_{\omega}$,
$A_{\mathcal{F}}(\omega)\subset A_{\mathcal{C}}(\omega)$. Thus the
result of this theorem follows. \qquad \
$\Box$\\

\section{Corollaries}
\ \ \;In this section useful corollaries following the main theorem
is
presented.\\
 \textbf{Corollary 4.1} If $\mathcal{C}$ is a
semi-algebra on $\Omega$, and $\mathcal{F}=\sigma(\mathcal{C})$.
Then $\forall$ $\omega\in\Omega$, $$A_{\mathcal{C}}(\omega)=
A_{\mathcal{F}}(\omega). $$ \textbf{Proof:} For any
$A\in\mathcal{C}$, one has
$$A^{c}=\Omega/A\in\mathcal{C}_{\sum f}\subset\mathcal{C}_{\sigma}
\subset\kappa(\mathcal{C}_{\sigma})=\kappa(\mathcal{C}). $$
Hence the result follows. \ \ \ \ $\Box$\\\\
\textbf{Corollary 4.2} If $\mathcal{C}$ is a semi-ring,
$\Omega\in\mathcal{C}_{\sigma}$, $\mathcal{F}=\sigma(\mathcal{C})$.
Then $\forall$ $\omega\in\Omega$,
$$A_{\mathcal{C}}(\omega)= A_{\mathcal{F}}(\omega). $$
\textbf{Proof:} There exists a sequence $A_{n}\in\mathcal{C}$ such
that $\Omega=\bigcup\limits_{n=1}^{\infty}A_{n}$. Then
$$A^{c}=\bigcup\limits_{n=1}^{\infty}(A_{n}/A), $$ by noting
$A_{n}/A\in\mathcal{C}_{\sum f}\subset\mathcal{C}_{\sigma}$. Hence
$A^{c}\in\mathcal{C}_{\sigma}\subset\kappa(\mathcal{C})$.\ \ \ \
$\Box$\\\\
\textbf{Corollary 4.3} If $\kappa(\mathcal{C})=\sigma(\mathcal{C})$,
then
 $\forall$ $\omega\in\Omega$,
$A_{\mathcal{C}}(\omega)= A_{\mathcal{F}}(\omega)$.\\
\textbf{Proof:} We show the equivalence between
$\kappa(\mathcal{C})=\sigma(\mathcal{C})$ and $\forall$
$A\in\mathcal{C}$ $\Rightarrow$ $A^{c}\in \kappa(\mathcal{C})$. If
$\kappa(\mathcal{C})=\sigma(\mathcal{C})$,  obviously, there are
$\forall$ $A\in\mathcal{C}$ $\Rightarrow$ $A^{c}\in
\kappa(\mathcal{C})$. For the inverse direction,  the collection of
set $\mathcal{G} =\{B\in\kappa(\mathcal{C})|$
$B^{c}\in\kappa(\mathcal{C})\}$, which is
 closed under countable intersection, countable union and
 complement, contains $\mathcal{C}$. Hence, $\mathcal{G}$ is a
  $\sigma$ algebra and $\mathcal{G}=\kappa
(\mathcal{C})=\sigma(\mathcal{C})$.\ \ \ \
$\Box$.\\\\
In the following corollaries we suppose that $\mathcal{F}$ is
separable($\mathcal{F}$ can be generated by a countable subset), so
they can be directly applied to the comparison among atoms in
Blackwell space\cite{DM78}.\\\\
\textbf{Corollary 4.4} Suppose $\mathcal{C}$ is a countable
semi-ring and $\mathcal{F}=\sigma(\mathcal{C})$(Obviously,
$\mathcal{F}$ is separable). Then $\forall$ $\omega\in\Omega$,
$A_{\mathcal{C}}(\omega)=
A_{\mathcal{F}}(\omega)$ if and only if $\Omega\in\mathcal{C}_{\sigma}$.\\
\textbf{Proof:} From Corollary 4.2, we know we only have to check if
$\forall$ $\omega\in\Omega$ $A_{\mathcal{C}}(\omega)=
A_{\mathcal{F}}(\omega)$ $\Rightarrow$
$\Omega\in\mathcal{C}_{\sigma}$. $\forall$ $\omega\in\Omega$
$A_{\mathcal{C}}(\omega)= A_{\mathcal{F}}(\omega)$, then $\exists$
$B\in\mathcal{C}$ such that $\omega\in B$. Hence,
$\Omega=\bigcup\limits_{B\in\mathcal{C}}B$. Note $\mathcal{C}$ is
countable, then $\Omega=\bigcup\limits_{B\in\mathcal{C}}B\in
\mathcal{C}_{\sigma}$.\ \ \ \
$\Box$ \\\\
\textbf{Corollary 4.5} Let $\mathcal{C}$ be a collection of sets on
$\Omega$, $\mathcal{F}=\sigma(\mathcal{C})$. If $\mathcal{F}$ has
countable atoms and $\mathcal{C}$ is countable. Then $\forall$
$\omega\in\Omega$, $A_{\mathcal{F}}(\omega)=A_{\mathcal{C}}(\omega)$
if and only if  $\mathcal{F}=\kappa(\mathcal{C})$.\\
\textbf{Proof:} If $\forall$ $\omega\in\Omega$,
$A_{\mathcal{F}}(\omega)=A_{\mathcal{C}}(\omega)$.  $\forall$
$A\in\mathcal{C}$, $A^{c}=\bigcup\limits_{\omega\in
A^{c}}A_{\mathcal{F}}(\omega)=\bigcup\limits_{\omega\in
A^{c}}A_{\mathcal{C}}(\omega)$. Since $\mathcal{C}$ is countable,
$A_{\mathcal{C}}(\omega)=\bigcap\limits_{B\in\mathcal{C}_{\omega}}B\in\kappa(\mathcal{C})$.
Since the atoms of $\mathcal{F}$ is countable,
$\bigcup\limits_{\omega\in
A^{c}}A_{\mathcal{C}}(\omega)=\bigcup\limits_{\omega\in
A^{c}}\bigcap\limits_{B\in\mathcal{C}_{\omega}}B\in\kappa(\mathcal{C})$,
indicating $A^{c}\in\kappa(\mathcal{C})$. From the proof of
Corollary 4.3, $A^{c}\in\kappa(\mathcal{C})$ ($\forall$
$A\in\mathcal{C}$)implies $\mathcal{F}=\kappa(\mathcal{C})$. The
inverse that $\mathcal{F}=\kappa(\mathcal{C})$ implies $\forall$
$\omega\in\Omega$ $A_{\mathcal{F}}(\omega)=A_{\mathcal{C}}(\omega)$
is trivial if we note that $\mathcal{F}=\kappa(\mathcal{C})$ implies
$\forall$ $A\in\mathcal{C}$, $A^{c}\in\kappa(\mathcal{C})$.\ \ \ \
$\Box$\\

In Corollary 4.1 and 4.2, we do not really use the property of
semi-ring or semi-algebra, which is closed under finite
intersection. Besides the condition  $A\cap
B^{c}\in\mathcal{C}_{\sum f}$ can be replaced by $A\cap B^{c}\in
\mathcal{C}_{\sigma\delta}$.

\section{Discussion and conclusion}

\ \ \;First consider the following two examples.\\
\textbf{Example 5.1} Let $\Omega=R$, $\mathcal{C}=\{x|$ $x\in R\}$
and $\mathcal{F}=\sigma(\mathcal{C})$. Obviously, $\forall$ $x\in R$
$A_{\mathcal{C}}(x)={x}=A_{\mathcal{F}}(x)$, and $\mathcal{F}$ is
Hausdoff(the atoms of $\mathcal{F}$ are the points of $\Omega$).
It's easy to check $\kappa(\mathcal{C})\subset\{A\subset R|$ $A$ is
countable $\}$. However, $R/\{0\}\in \mathcal{F}$ is not in
$\kappa(\mathcal{C})$. This shows our condition is not a
necessary one.\\\\
\textbf{Example 5.2} Let $\Omega=[0,1]$, $\mathcal{C}=\{[a,b)\subset
[0,1)|$ $a<b\}\cup\{\emptyset\}$ and
$\mathcal{F}=\sigma(\mathcal{C})$. $\mathcal{C}$ is a semi-ring on
[0,1]. $A_{\mathcal{C}}(1)=\emptyset$, while
$A_{\mathcal{F}}(1)=\{1\}$. This shows the condition that the
generator is a semi-ring is not sufficient for
$A_{\mathcal{C}}(\omega)=A_{\mathcal{F}}(\omega)$($\forall$
$\omega\in\Omega$).\\

The examples show that our condition may not be the best one but
almost necessary. Comparing our condition with semi-ring (see
Corollary 4.2), we only add $\Omega\in\mathcal{C}_{\sigma}$ to
obtain the desired result, and the condition of Corollary 4.2 is
stronger than that in our main theorem. Therefore our condition has
already been a very weak one. On the other hand since
$m(\mathcal{C})\subset
m(\mathcal{C}_{\sigma})\subset\kappa(\mathcal{C}_{\sigma})=\kappa(\mathcal{C})
\subset\sigma(\mathcal{C})$ and $m(\mathcal{C}_{\sigma})\subset
\lambda(\mathcal{C}_{\sigma})\subset\sigma(\mathcal{C})$. The
trivial case $A^{c}\in\sigma(\mathcal{C})$ contributes nothing if
letting it replace $A^{c}\in \kappa(\mathcal{C})$ since it is
impossible to conclude
$A_{\mathcal{C}}(\omega)=A_{\mathcal{F}}(\omega)$($\forall$
$\omega\in\Omega$) without any restriction on $\mathcal{C}$ (Example
5.1 can be viewed as a special counterexample). From relations among
$m(\mathcal{C})$, $m(\mathcal{C}_{\sigma})$,
$\kappa(\mathcal{C}_{\sigma})$, $\kappa(\mathcal{C}) $, $
\lambda(\mathcal{C}_{\sigma})$, $\sigma(\mathcal{C})$, we know
$\kappa(\mathcal{C})$ is already a very large set and nearly as
large as $\lambda(\mathcal{C}_{\sigma})$. Finally, rewiewing the
proof of the theorem, one can find the key of the proof lies in the
property of $\mathcal{G}$ in Lemma 3.1. Generally, $\mathcal{G}$ is
at most a $\kappa$ class and could not be a $\lambda$ class. Hence
the improvement of our condition from the theoretic
perspective is almost impossible.\\

 {\bf Acknowledgements.} The author is very grateful to Professor Yves Le Jan
for helpful discussions and bringing the reference \cite{DM78} to
his attention.

\end{document}